\documentclass[10pt,a4paper, xcolor]{article}
\usepackage{amsfonts}
\usepackage{amssymb,amsthm, amsmath,graphicx,bm,amsthm}
\usepackage{color,enumerate,lineno}
\usepackage{url,hyperref}%,showlabels}

\usepackage{amsfonts}
\usepackage{amssymb,amsthm, amsmath,graphicx,bm,amsthm,mathtools}
\usepackage{color,enumerate}
\usepackage{url,lineno}

\usepackage[latin1]{inputenc}

\usepackage[all]{xy}

\usepackage{algorithm}
\usepackage[noend]{algpseudocode}

\newtheorem{theorem}{Theorem}

\newtheorem{proposition}[theorem]{Proposition}
\newtheorem{corollary}[theorem]{Corollary}
\newtheorem{lemma}[theorem]{Lemma}

\theoremstyle{remark}

\newtheorem{example}[theorem]{Example}

\def\gcd{{\mathrm gcd}}

\def\N{\mathbb{N}}

\def\Z{\mathbb{Z}}

\def\Ap{\mathrm{Ap}}
\def\int{\mathrm{int}}

\DeclarePairedDelimiter\ceil{\lceil}{\rceil}
\DeclarePairedDelimiter\floor{\lfloor}{\rfloor}

\title{Frobenius number and minimum genus of numerical semigroups with fixed multiplicity and embedding dimension}

\author{
J. I. Garc\'{\i}a-Garc\'{\i}a
\footnote{
	Dpto. de Matem\'aticas/INDESS (Instituto Universitario para el Desarrollo Social Sostenible).
	Universidad de C\'adiz, E-11510 Puerto Real  (C\'{a}diz, Spain).
	E-mail: ignacio.garcia@uca.es.
	Partially supported by MTM2014-55367-P and by Junta de Andaluc\'{\i}a group FQM-366.
}
\and
D. Mar\'{\i}n-Arag\'{o}n
\footnote{
	E-mail: daniel.marinaragon@alum.uca.es.
	Partially supported by Junta de Andaluc\'{\i}a group FQM-366.
}
\and
M. A. Moreno-Fr\'{\i}as
\footnote{
	Dpto. Matem\'aticas,
	Universidad de C\'adiz, E-11510 Puerto Real  (C\'{a}diz, Spain).
	E-mail: mariangeles.moreno@uca.es.
	Partially supported by MTM2014-55367-P and by Junta de Andaluc\'{\i}a group FQM-298.
}
\and
J. C. Rosales
\footnote{
	Dpto. de \'Algebra, Universidad de Granada.
	Partially supported by MTM2014-55367-P and by Junta de Andaluc\'{\i}a group FQM-343.
	E-mail: jrosales@ugr.es.
}
\and
A. Vigneron-Tenorio
\footnote{
	Dpto. de Matem\'aticas/INDESS (Instituto Universitario para el Desarrollo Social Sostenible).
	Universidad de C\'adiz,
 	E-11406 Jerez de la Frontera (C\'{a}diz, Spain).
 	E-mail: alberto.vigneron@uca.es.
 	Partially supported by MTM2015-65764-C3-1-P (MINECO/FEDER, UE) and
	Junta de Andaluc\'{\i}a group FQM-366. 	
}
}

\date{}

\begin{document}

\maketitle

\begin{abstract}
Fixed two
%natural numbers
positive integers
$m$ and $e,$ some algorithms for computing the minimal Frobenius number and minimal genus of the set of numerical semigroups with multiplicity $m$ and embedding dimension $e$ are provided. Besides, the semigroups where these minimal values are achieved are computed too.


    \smallskip
    {\small \emph{MSC-class:} 20M14 (Primary),  20M05 (Secondary).}
\end{abstract}

\section{Introduction}
Given a subset $S$ of the set of non-negative integers denoted by $\N,$ $S$ is called a numerical semigroup if $0\in S,$ $\N\setminus S$ is finite and $S$ is closed by sum.
%
%Let $\N=\{0,1,2,\ldots\}$ be the set of non-negative integers. A numerical semigroup is a subset $S$ of $\N$ which is closed by sum, $0\in S$ and $\N\setminus S$ is finite.

If $S$ is a numerical semigroup, it is well-known there exists a unique inclusion-wise minimal finite subset $A$ of $\N$ such that $S=\{\lambda_1 a_1+\cdots+\lambda_n a_n \mid n\in\N, ~ a_1,\ldots,a_n\in A$ and $\lambda_1,\ldots,\lambda_n\in\N\}$ (see Theorem 2.7 from \cite{libro}). The set $A$ is called the minimal system of generators of $S.$ In general, given a finite set $A\subset \N,$ the monoid generated by $A$, $\langle A \rangle= \{\lambda_1 a_1+\cdots+\lambda_n a_n \mid n\in\N, ~ a_1,\ldots,a_n\in A$ and $\lambda_1,\ldots,\lambda_n\in\N\},$ is a numerical semigroup if and only if $\gcd(A)=1$ (see \cite[Lemma 2.1]{libro}). We denote by $msg(S)$ the minimal system of generators of $S$; its cardinality is called the embedding dimension of $S$ and it is denoted by $e(S)$.

%Given a non-empty subset $A$ of $\N$ we denote by $\langle A \rangle$ to the submonoid of $(\N,+)$ generated by $A$, that is, $\langle A \rangle = \{\lambda_1 a_1+\cdots+\lambda_n a_n \mid n\in\N, ~ a_1,\ldots,a_n\in A$ and $\lambda_1,\ldots,\lambda_n\in\N\}$. It is well known  that $\langle A \rangle$ is a numerical semigroup if and only if $\gcd(A)=1$ (see, for example, Lemma 2.1 from \cite{libro}).
%
%If $S$ is a numerical semigroup and $S=\langle A \rangle$ then we say that $A$ is a system of generators of $S$. Moreover, if $S\neq\langle B \rangle$ for every $B\subsetneq A$ then we say that $A$ is a minimal system of generators of $S$. In Theorem 2.7 from \cite{libro} it is shown that every numerical semigroup has a unique minimal system of generators and this system is finite. We denote by $msg(S)$ the minimal system of generators of $S$; its cardinality is called the embedding dimension of $S$ and it is denoted by $e(S)$.

Other important invariants of a numerical semigroup $S$ are the following:
\begin{itemize}
    \item The multiplicity of $S$ (denoted by $m(S)$) is the least positive integer that belongs to $S$.
    \item The Frobenius number of $S$ (denoted by $F(S)$) is the greatest integer which does not belong to $S$.
    \item The genus of $S$ (denoted by $g(S)$) is the cardinality of $\N\setminus S$.
\end{itemize}

For those readers who are not familiarized with Numerical Semigroup Theory, the terminology used (genus, multiplicity, embedding dimension,\dots) may result bizarre in the scope of semigroups.
In the literature, one finds many manuscripts devoted
to the study of analytically unramified one dimensional local domains via their
value semigroups, which turn out to be numerical semigroups.
All these invariants (see \cite{barucci}) have interpretations in this theory, and hence its names.

The relation of these invariantes is of great interest. Indeed, in \cite{wilf} the following conjeture appears: if $S$ is a numerical semigroup then $e(S)g(s)\leq (e(S)-1)(F(S)+1)$. Nowadays, this conjecture is still open and it is one of the most important problems in Numerical Semigroup Theory.

If $m$ and $e$ are positive integers we give the following notations:
$\mathcal{L}(m,e)=\{S \mid   S \mbox{ is a numerical semigroup, }  m(S)=m, ~ e(S)=e \}$,
$g(m,e)=\min\{g(S) \mid S\in\mathcal{L}(m,e) \}$ and $F(m,e)=\min\{F(S) \mid S\in\mathcal{L}(m,e) \}$. In this work we are interested in giving algorithms for computing $g(m,e)$, $F(m,e)$, $\{S\in\mathcal{L}(m,e) \mid g(S)=g(m,e) \}$ and $\{S\in\mathcal{L}(m,e) \mid F(S)=F(m,e) \}$.

The results of this work are illustrated with several examples. To this aim, we have used the library \texttt{FrobeniusNumberAndGenus} developed by the authors in Mathematica (\cite{mathematica}). The library \texttt{FrobeniusNumberAndGenus} is freely available online at \cite{PROGRAMA}

This work is organized as follows. Section \ref{pseudo} contains some known results on Frobenius pseudo-varieties which allow us to construct the tree of numerical semigroup with fixed multiplicity. In Section \ref{mingenus} and Section \ref{minFrob}, the minimal genus and minimal Frobenius number of the set of numerical semigroups with fixed multiplicity and embedding dimension are studied, giving some algorithm for computing them and obtaining the semigroups with these minimal values.

\section{Frobenius pseudo-variety of numerical semigroups with a fixed multiplicity}\label{pseudo}

According to the notation of \cite{pseudo-variedades}, a Frobenius pseudo-variety is a non-empty family $\mathcal{P}$ of numerical semigroups which verifies the following conditions:
\begin{enumerate}
    \item $\mathcal{P}$ has a maximum (according to the inclusion order).
    \item If $\{S,T\}\subseteq\mathcal{P}$ then $S\cap T\in\mathcal{P}$.
    \item If $S\in\mathcal{P}$ and $S\neq\max(\mathcal{P})$ then $S\cup\{F(S)\}\in\mathcal{P}$.
\end{enumerate}

A graph $G$ is a pair $(V,E)$ where $V$ is a set (vertices of $G$) and $E$ is a subset of $\{(u,v)\in V\times V \mid u\neq v\}$ (edges of $G$). A path of length $n$ connecting the vertices $x$ and $y$ is a sequence of different edges of the form $(v_0,v_1), (v_1,v_2),\ldots,(v_{n-1},v_n)$ such that $v_0=x$ and $v_n=y$.

A graph $G$ is a tree if there exists a vertex $r$ (known as the root) such that for any other vertex $x$ of $G$ there exists a unique path connecting $x$ and $r$. If $(x,y)$ is an edge of a tree, we say that $x$ is a son of $y$.

If $\mathcal{P}$ is a Frobenius pseudo-variety we define the graph $G(\mathcal{P})$ as follows: $\mathcal{P}$ is its set of vertices and $(S,T)\in\mathcal{P}\times\mathcal{P}$ is an edge if $T=S\cup\{F(S)\}$.

The following result is a direct consequence of Lemma 12 and Theorem 3 of \cite{pseudo-variedades}.
% As a consequence of of Lemma 12 and Theorem 3 of \cite{pseudo-variedades}, we obtain the following result.

\begin{proposition}\label{p1}
If $\mathcal{P}$ is a Frobenius pseudo-variety, then $G(\mathcal{P})$ is a tree with root $\max(\mathcal{P})$. Moreover, the set of sons of a vertex $S\in\mathcal{P}$ is $\{S\setminus\{x\}\in\mathcal{P} \mid x\in msg(S),\ x>F(S) \}$.
\end{proposition}

Let $m$ be a positive integer. We denote by $\mathcal{L}(m)$ the set \[\{S \mid S \mbox{ is a numerical semigroup with } m(S)=m \}.\] Clearly $\mathcal{L}(m)$ is a Frobenius pseudo-variety and $\max(\mathcal{L}(m))=\{0,m,\rightarrow\}=\langle m,m+1,\ldots,2m-1 \rangle$. So, as a consequence of Proposition \ref{p1} we have the following result which is fundamental in this work.

\begin{theorem}\label{t2}
The graph $G(\mathcal{L}(m))$ is a tree rooted in $\langle m,m+1,\ldots,2m-1 \rangle$. Moreover, the set formed by the sons of a vertex $S\in\mathcal{L}(m)$ is $\{S\setminus\{x\} \mid x\in msg(S),\ x>F(S) \mbox{ and } x\neq m \}$.
\end{theorem}

The previous theorem allows us to build recursively $\mathcal{L}(m)$ from its root and adding to the computed vertices its sons.
%(which are characterized in Theorem \ref{t2}).
We illustrate it with an example.

\begin{example}\label{e3}
We show some levels of tree $G(\mathcal{L}(4))$ giving its vertices and edges, and the minimal removed generators for obtaining the sons.

%\hspace{-.5cm}
{\footnotesize
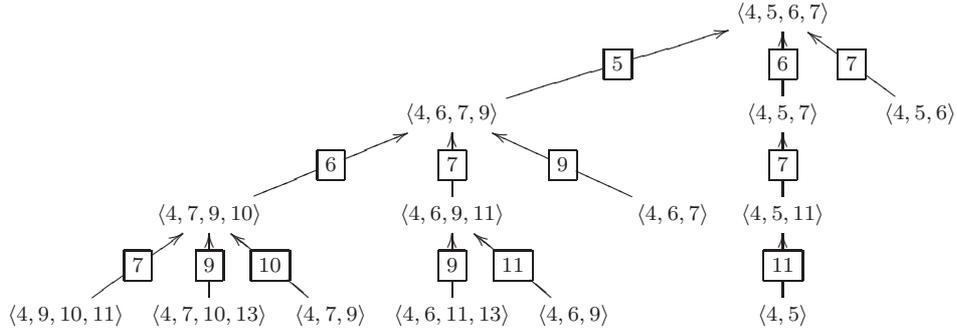
\xymatrix@C=0.5em{
 & & & & & & \langle 4,5,6,7\rangle \ar@{<-}[rrd]|*+[F]{7}\ar@{<-}[llld]|*+[F]{5}\ar@{<-}[d]|*+[F]{6}  & & & & &       \\
 & & & \langle 4,6,7,9 \rangle \ar@{<-}[rrd]|*+[F]{9}\ar@{<-}[lld]|*+[F]{6}\ar@{<-}[d]|*+[F]{7} & & & \langle 4,5,7\rangle \ar@{<-}[d]|*+[F]{7} & & \langle 4,5,6\rangle \\
 & \langle 4,7,9,10\rangle \ar@{<-}[rd]|*+[F]{10}\ar@{<-}[ld]|*+[F]{7}\ar@{<-}[d]|*+[F]{9}  & & \langle 4,6,9,11 \rangle \ar@{<-}[d]|*+[F]{9}\ar@{<-}[rd]|*+[F]{11} & & \langle 4,6,7 \rangle & \langle 4,5,11 \rangle \ar@{<-}[d]|*+[F]{11} & &            \\
\langle 4,9,10,11 \rangle & \langle 4,7,10,13 \rangle & \langle 4,7,9 \rangle  & \langle 4,6,11,13 \rangle & \langle 4,6,9 \rangle &  & \langle 4, 5\rangle & &     \\
}
}
\end{example}

If $G$ is a tree with root $r$, the level of a vertex $x$ is the length of the only path which connect $x$ and $r$. The height of a tree is the value of its maximum level. If $k\in\N$, we denote by $N(k,G)=\{v\in G \mid v \mbox{ has level  }k \}$. So, for Example \ref{e3} we have:

$$N(0,\mathcal{L}(4))=\{\langle 4,5,6,7 \rangle\}.$$
$$N(1,\mathcal{L}(4))=\{\langle 4,6,7,9 \rangle,\langle 4,5,7 \rangle,\langle 4,5,6 \rangle\}.$$
$$N(2,\mathcal{L}(4))=\{\langle 4,7,9,10 \rangle,\langle 4,6,9,11 \rangle,\langle 4,6,7 \rangle,\langle 4,5,11 \rangle\}.$$
$$N(3,\mathcal{L}(4))=\{\langle 4,9,10,11 \rangle,\langle 4,7,10,13 \rangle,\langle 4,7,9 \rangle,\langle 4,6,11,13 \rangle,\langle 4,6,9 \rangle,\langle 4,5 \rangle\}.$$

\section{Elements of $\mathcal{L}(m,e)$ with minimum genus}\label{mingenus}

Our aim in this section is to give an algorithm that allows us to compute $g(m,e)$ and $\{S \mid S\in\mathcal{L}(m,e) \mbox{ and } g(S)=g(m,e) \}$. The following result is a consequence of Theorem \ref{t2}.

\begin{proposition}
If $m$ is a positive integer and $(S,T)$ an edge of $G(\mathcal{L}(m))$, then $g(S)=g(T)+1$.
\end{proposition}

As a direct consequence of the previous proposition we have the following result.

\begin{corollary}
Fixed $m,e\in \N,$ if $P=\min\{k\in\N \mid N(k,G(\mathcal{L}(m)))\cap\mathcal{L}(m,e)\neq\emptyset\}$ then $\{S\in\mathcal{L}(m,e) \mid g(S)=g(m,e)\}=N(P,G(\mathcal{L}(m)))\cap\mathcal{L}(m,e)$. Moreover, $g(m,e)=m-1+P$.
\end{corollary}

Note that if $S$ is a numerical semigroup and $e(S)=1$ then $S=\N$. Note also that by Proposition 2.10 of \cite{libro}, if $S$ is a numerical semigroup then $e(S)\leq m(S)$. Moreover, it is clear that if $m\geq e\geq 2$ then $\langle m, m+1,\ldots,m+e-1\rangle\in\mathcal{L}(m,e)$. In this way, we have the following result.

\begin{proposition}
Let $m$ and $e$ be positive integers.
    \begin{enumerate}
    \item If $m<e$ then $\mathcal{L}(m,e)=\emptyset$.
    \item If $e=1$ and $\mathcal{L}(m,e)\neq\emptyset$ then $m=1$ and $\mathcal{L}(m,e)=\{\N\}$.
    \item If $m\geq e\geq 2$ then $\mathcal{L}(m,e)\neq\emptyset$.
    \end{enumerate}
\end{proposition}

Now, we give an algorithm to compute $g(m,e)$ and $\{S\in\mathcal{L}(m,e) \mid g(S)=g(m,e) \}$.

\begin{algorithm}
\caption{Sketch of the algorithm to compute $g(m,e)$ and the set of semigroups with a fixed multiplicity and embedding dimension such that its genus is $g(m,e)$.}\label{a7}
\textbf{INPUT:} $m$ and $e$ positive integers such that $m\geq e\geq 2$.\\
\textbf{OUTPUT:} $g(m,e)$ and $\{S \mid S\in\mathcal{L}(m,e) \mbox{ and } g(S)=g(m,e) \}$.
\begin{algorithmic}[1]
	\State Set $k=0$ and $A=\{\langle m,m+1,\ldots,2m-1\rangle\}$
	\While{True}
		\If {$A\cap\mathcal{L}(m,e)\neq\emptyset$}
			\State \Return{$m-1+k$ and $A\cap\mathcal{L}(m,e)$}
		\EndIf
		\For{$S\in A$}
			\State $C(S)=\{T \mid T \mbox{ is a son of } S\}$
		\EndFor
		\State $\displaystyle{A=\bigcup_{S\in A}C(S)}$, $k=k+1$.
	\EndWhile
\end{algorithmic}
%\begin{enumerate}[Step 1.]
%\item Set $k=0$ and $A=\{m,m+1,\ldots,2m-1\}$.
%\item If $A\cap\mathcal{L}(m,e)\neq\emptyset$ then return $m-1+k$ and $A\cap\mathcal{L}(m,e)$.
%\item For every $S\in A$, $C(S)=\{T \mid T \mbox{ is a son of } S\}$.
%\item $\displaystyle{A=\bigcup_{S\in A}C(S)}$, $k=k+1$ and go to 2.
%\end{enumerate}
\end{algorithm}

We illustrate now the previous algorithm with an example.

\begin{example}
We compute $g(5,3)$ and $\{S\in\mathcal{L}(5,3) \mid g(S)=g(5,3) \}$ using Algorithm \ref{a7}.
    \begin{itemize}
    \item $k=0$ and $A=\{\langle 5,6,7,8,9 \rangle\}$.
    \item $k=1$ and $A=\{\langle 5,7,8,9,11 \rangle,\langle 5,6,8,9 \rangle,\langle 5,6,7,9 \rangle,\langle 5,6,7,8 \rangle\}$.
    \item $k=2$ and
    \begin{multline*}
    A=\{\langle 5,8,9,11,12 \rangle,\langle 5,7,9,11,13 \rangle,\langle 5,7,8,11 \rangle,\\
    \langle 5,7,8,9 \rangle, \langle 5,6,9,13 \rangle,\langle 5,6,8 \rangle,\langle 5,6,7 \rangle\}.
    \end{multline*}
    \end{itemize}
It returns $g(5,3)=6$ and $\{S\in\mathcal{L}(5,3) \mid g(S)=6 \}=\{\langle 5,6,8 \rangle,\langle 5,6,7 \rangle\}$.
In the package \texttt{FrobeniusNumberAndGenus}, we can run the command \texttt{Algorithm1[5,3]} to obtain this result.
\end{example}

If $S$ is a numerical semigroup and $n\in S\setminus\{0\}$, then the Ap\'{e}ry set of $n$ in $S$ is $\Ap(S,n)=\{s\in S \mid s-n\not\in S\}$ (see \cite{apery}). Lemma 2.4 of \cite{libro} shows that $\Ap(S,n)=\{w(0)=0,w(1),\ldots,w(n-1)\}$ where $w(i)$ is the least element of $S$ congruent with $i$ modulus $n$. Note that $w(i)=k_in+i$ for some $k_i\in\N$ and $km+i\in S$ if and only if $k\geq k_i$. Therefore, we have the following result.

\begin{lemma}\label{l9}
Let $S$ be a numerical semigroup, $n\in S\setminus\{0\}$ and $\Ap(S,n)=\{0,k_1n+1,\ldots,k_{n-1}n+n-1\}$. Then $g(S)=k_1+\cdots+k_{n-1}$.
\end{lemma}

Next result can be easily deduced from Corollary 4 of \cite{intervalos}.

\begin{lemma}\label{l10}
Let $m$ and $e$ integers such that $m\geq e\geq 2$, $S=\langle m,m+1,\ldots,m+e-1 \rangle$ and $m-1=q(e-1)+r$ with $q,r\in\N$ and $r\leq e-2$. Then $\Ap(S,m)=\{0,m+1,\ldots,m+e-1,2m+(e-1)+1,\ldots,2m+2(e-1),\ldots,qm+(q-1)(e-1)+1,\ldots,qm+q(e-1),(q+1)m+q(e-1)+1,\ldots,(q+1)m+q(e-1)+r\}$.
\end{lemma}

If $a,b\in\N$ and $b\neq 0$ we denote by $a\mod b$ the remainder of dividing $a$ by $b$. If $q$ is a rational number we denote by $\floor{q}=\max\{z\in\Z \mid z\leq q\}$. Note that $a=\floor{\frac a b} b+(a\mod b)$. From Lemma \ref{l9} and Lemma \ref{l10} we have the following result.

\begin{proposition}
Let $m$ and $e$ be integers such that $m\geq e\geq 2$ and $S=\langle m,m+1,\ldots,m+e-1 \rangle$. Then, $$g(S)=\left(\left\lfloor{\frac{m-1}{e-1}+1}\right\rfloor\right)\left(\frac{\left\lfloor{\frac{m-1}{e-1}}\right\rfloor(e-1)}{2}+(m-1)\mod (e-1)\right).$$
\end{proposition}

Clearly $\langle m,m+1,\ldots,m+e-1 \rangle \in\mathcal{L}(m,e)$ and therefore we have the following result.

\begin{corollary}
If $m$ and $e$ are integers such that $m\geq e\geq 2$ then $$g(m,e)\leq\left(\left\lfloor{\frac{m-1}{e-1}+1}\right\rfloor\right)\left(\frac{\left\lfloor{\frac{m-1}{e-1}}\right\rfloor(e-1)}{2}+(m-1)\mod (e-1)\right).$$
\end{corollary}

For many examples the equality holds. However, there are some cases where the semigroup $\langle m,m+1,\ldots,m+e-1 \rangle$ does not have minimum genus in the set $\mathcal{L}(m,e)$ as we show in the next example.

\begin{example}
$S=\langle 8,9,10 \rangle$ is a numerical semigroup and $g(S)=16$. $\bar{S}=\langle 8,9,11 \rangle$ is a numerical semigroup and $g(\bar{S})=14$. Therefore, in this case $g(\langle 8,9,10 \rangle)\neq g(8,3)$.
\end{example}

Using the notation of \cite{cadiz}, a numerical semigroup is packed if $msg(S)\subseteq\{m(S),m(S)+1,\ldots,2m(S)-1\}$. The set of all packed numerical semigroups with multiplicity $m$ and embedding dimension $e$ is denote by $\mathcal{C}(m,e)$. The following result is obtained from \cite{cadiz}.

\begin{proposition}\label{p14}
If $S\in\mathcal{L}(m,e)$ then $\bar{S}=\langle \{m\}+\{x\mod m \mid x\in msg(S)\} \rangle\in\mathcal{C}(m,e)$ and $g(\bar{S})\leq g(S)$. Moreover, if $S\not\in\mathcal{C}(m,e)$ then $g(\bar{S})<g(S)$.
\end{proposition}

We illustrate the previous proposition with an example.

\begin{example}
If $S=\langle 5,11,17 \rangle\in\mathcal{L}(5,3)$ then $\bar{S}=\langle \{5\}+\{0,1,2\} \rangle=\langle 5,6,7 \rangle\in\mathcal{C}(5,3)$. Therefore, $g(\bar{S})\leq g(S)$. Moreover, $S\not\in\mathcal{C}(5,3)$, so $g(\bar{S})<g(S)$.
\end{example}

Next result is a consequence of Proposition \ref{p14}.

\begin{corollary}
Let $m$ and $e$ be integers such that $m\geq e \geq 2$. Then
    \begin{enumerate}
        \item $g(m,e)=\min\{g(s) \mid S\in\mathcal{C}(m,e)\}$.
        \item $\{S\in\mathcal{L}(m,e) \mid g(S)=g(m,e) \}=\{S\in\mathcal{C}(m,e) \mid g(S)=g(m,e) \}$.
    \end{enumerate}
\end{corollary}

Note that $\mathcal{C}(m,e)$ is finite and therefore the previous corollary give us another algorithm for computing $g(m,e)$ and $\{S\in\mathcal{L} \mid g(S)=g(m,e) \}$. We give more details about this method. Following result is Proposition 4 of \cite{cadiz}.

\begin{proposition}\label{p17}
Let $m$ and $e$ integers such that $m\geq e\geq 2$ and $A$ a subset of $\{1,\ldots,m-1\}$ with cardinality $e-1$ such that $gcd(A\cup\{m\})=1$. Then $S=\langle \{m\}+(A\cup\{0\}) \rangle\in\mathcal{C}(m,e)$. Moreover, every element of $\mathcal{C}$ has this form.
\end{proposition}

We illustrate the previous proposition with an example.

\begin{example}
We are going to compute $\mathcal{C}(6,3)$. We start computing $\{A\subseteq\{1,2,3,4,5\} \mid \#A=2 \mbox{ and } \gcd(A\cup\{6\})=1 \}$. This set is equal to
\[
\{ \{1,2\},\{1,3\},\{1,4\},\{1,5\},\{2,3\},\{2,5\},\{3,4\}, \{3,5\},\{4,5\}\}.
\]

Therefore,
\begin{multline*}
\mathcal{C}(6,3)=\{\langle 6,7,8\rangle,\langle 6,7,9\rangle,\langle 6,7,10\rangle,\langle 6,7,11\rangle,\\
\langle 6,8,9\rangle,\langle 6,8,11\rangle,\langle 6,9,10\rangle,\langle 6,9,11\rangle,\langle 6,10,11\rangle\}.
\end{multline*}
A simple computation shows us $g(\langle 6,7,8 \rangle)=9$, $g(\langle 6,7,9 \rangle)=9$, $g(\langle 6,7,10 \rangle)=9$, $g(\langle 6,7,11 \rangle)=10$, $g(\langle 6,8,9 \rangle)=10$, $g(\langle 6,8,11 \rangle)=11$, $g(\langle 6,9,10 \rangle)=12$, $g(\langle 6,9,11 \rangle)=13$ and $g(\langle 6,10,11 \rangle)=13$.\\
Therefore, $g(6,3)=9$ and the set $\{S\in\mathcal{L}(6,3) \mid g(S)=9 \}$ is equal to $\{\langle 6,7,8 \rangle,\langle 6,7,9 \rangle,\langle 6,7,10 \rangle\}$.
\end{example}

\section{Elements of $\mathcal{L}(m,e)$ with minimum Frobenius Number}\label{minFrob}

Our aim in this section is to obtain algorithmic methods for computing $F(m,e)$ and $\{S\in\mathcal{L}(m,e) \mid F(S)=F(m,e) \}$. Next result is a consequence of Theorem \ref{t2}.

\begin{proposition}
If $m$ is a positive integer and $(S,T)$ is an edge of $G(\mathcal{L}(m))$ then $F(T)<F(S)$.
\end{proposition}

The following result can be deduced from \cite{numerical}.

\begin{proposition}
If $m$ is a positive integer and $(S,T)$ is an edge of $G(\mathcal{L}(m))$ then $e(S)\leq e(T)$.
\end{proposition}

Clearly $F(m,m)=m-1$ and $\{S\in\mathcal{L}(m,m) \mid F(S)=m-1 \}=\{\langle m,m+1,\ldots,2m-1 \rangle\}$. It is well known (see \cite{sylvester} for example) that if $S=\langle n_1,n_2 \rangle$ is a numerical semigroup,
%of embedding dimension two
then $F(S)=n_1n_2-n_1-n_2$. Therefore, we obtain the following result.

\begin{proposition}
Let $m$ be an integer such that $m\geq 2$.
    \begin{enumerate}
        \item $F(m,m)=m-1$ and $\{S\in\mathcal{L}(m,m) \mid F(S)=m-1\}=\{\langle m,m+1,\ldots,2m-1\rangle\}$.
        \item $F(m,2)=m^2-m-1$ and $\{S\in\mathcal{L}(m,2) \mid F(S)=m^2-m-1\}=\{\langle m,m+1 \rangle\}$.
    \end{enumerate}
\end{proposition}

If $q$ is a rational number we denote by $\ceil{q}=\min\{z\in\Z \mid q\leq z \}$. Next result is deduced from \cite[Corollary \S 4.5]{intervalos}.

\begin{proposition}
If $m$ and $e$ are integers such that $m\geq e\geq 2$ then $F(\langle m,m+1,\ldots,m+e-1 \rangle)=\ceil{\frac{m-1}{e-1}}m-1$.
\end{proposition}

As a consequence of the previous proposition we get the following result.

\begin{corollary}
If $m$ and $e$ are integers such that $m\geq e\geq 2$ then $F(m,e)\leq \ceil{\frac{m-1}{e-1}}m-1$.
\end{corollary}

In the previous corollary, the equality holds plenty of times, but in some cases $F(\langle m,m+1,\dots,m+e-1 \rangle)\neq\min\{F(S) \mid S\in\mathcal{L}(m,e) \}$. For example, $F(\langle 4,5,6 \rangle)=7$ and $F(\langle 4,5,7 \rangle)=6$.

From the above results, we obtain the following algorithm where
the projections from the cartesian product $\mathcal{L}(m)\times \N$ are denoted by $\pi_1$ and $\pi_2$.

\begin{algorithm}
\caption{Sketch of the algorithm to compute $F(m,e)$ and the set of semigroup with a fixed multiplicity and embedding dimension such that its Frobenius number is $F(m,e)$.}\label{a24}
\textbf{INPUT:} $m$ and $e$ integers such that $m\geq e\geq 2$.\\
\textbf{OUTPUT:} $F(m,e)$ and $\{S\in\mathcal{L}(m,e) \mid F(S)=F(m,e) \}$.
\begin{algorithmic}[1]
	\State $A=\{\langle m,m+1,\ldots,2m-1 \rangle\}$, $I=\emptyset$ and $\alpha=\ceil{\frac{m-1}{e-1}}m-1$
	\While{True}
		\State $C=\{(S,F(S)) \mid S \mbox{ is son of some element of } A \mbox{ and } F(S)\leq\alpha \}$
		\State $K=\{S\in \pi_1(C) \mid e(S)\geq e \}$
		\If{$K=\emptyset$}
			\State \Return $F(m,e)=\pi_2(I)$ and $\{S\in\mathcal{L}(m,e) \mid F(S)=F(m,e) \}=\pi_1(I)$
		\EndIf
		\State $A=K$, $B=\{(S,F(S)) \mid S\in K \mbox{ and } e(S)=e \}$
		\State $\alpha=\min(\pi_2(B)\cup\{\alpha\})$, $I=\{(S,F(S))\in I\cup B \mid F(S)=\alpha \}$.
	\EndWhile
\end{algorithmic}

%    \begin{enumerate}[Step 1.]
%        \item $A=\{\langle m,m+1,\ldots,2m-1 \rangle\}$, $I=\emptyset$ and $\alpha=\ceil{\frac{m-1}{e-1}}m-1$.
%        \item $C=\{(S,F(S)) \mid S \mbox{ is son of some element of } A \mbox{ and } F(S)\leq\alpha \}$ and $K=\{S\in C \mbox{ such that } e(S)\geq e \}$.
%        \item If $K=\emptyset$ then return $F(m,e)=\pi_2(I)$ and $\{S\in\mathcal{L}(m,e) \mbox{ such that } F(S)=F(m,e) \}=\pi_1(I)$.
%        \item $A=K$, $B=\{(S,F(S)) \mid S\in K \mbox{ and } e(S))e \}$, $\alpha=\min(\pi_2(B)\cup\{\alpha\})$, $I=\{(S,F(S))\in I\cup B \mid F(S)=\alpha \}$ and go to Step 2.
%    \end{enumerate}
\end{algorithm}

We illustrate with an example how this algorithm works.

\begin{example}
We compute $F(4,3)$ and $\{S\in\mathcal{L}(4,3) \mid F(S)=F(4,3) \}$ using Algorithm \ref{a24}.
    \begin{itemize}
        \item $A=\{\langle 4,5,6,7 \rangle\}$, $I=\emptyset$ and $\alpha=\ceil{\frac{3}{2}}4-1=7$.
        \item $C=\{(\langle 4,6,7,9 \rangle,5),(\langle 4,5,7 \rangle,6),(\langle 4,5,6 \rangle,7)\}$ and\\
        $K=\{\langle 4,6,7,9 \rangle,\langle 4,5,7 \rangle,\langle 4,5,6 \rangle\}$.
        \item $A=\{\langle 4,6,7,9 \rangle,\langle 4,5,7 \rangle,\langle 4,5,6 \rangle\}$, $B=\{(\langle 4,5,7 \rangle,6),(\langle 4,5,6 \rangle,7)\}$, $\alpha=\min\{6,7\}=6$ and $I=\{(\langle 4,5,7 \rangle,6)\}$.
        \item $C=\{(\langle 4,7,9,10 \rangle,6)\}$ and $K=\{\langle 4,7,9,10 \rangle\}$.
        \item $A=\{\langle 4,7,9,10 \rangle\}$, $B=\emptyset$, $\alpha=6$ and $I=\{(\langle 4,5,7 \rangle,6)\}$.
        \item $C=\emptyset$ and $K=\emptyset$.
    \end{itemize}
Therefore, $F(4,3)=6$ and $\{S\in \mathcal{L}(4,3) \mid F(S)=6 \}=\{\langle 4,5,7 \rangle\}$. Using the Mathematica package \cite{PROGRAMA},  running the commands \texttt{MinFrob[4,3]} and \texttt{FrobeniusEmbeddingDimensionMultiplicity[6,3,4]}, we obtain $6$ and $\langle 4,5,7 \rangle$, respectively.
\end{example}

Our next goal is to give an alternative algorithm for computing $F(m,e)$ and $\{S\in\mathcal{L}(m,e) \mid F(S)=F(m,e) \}$. Next result is deduced from \cite{cadiz}.

\begin{proposition}
If $S\in\mathcal{L}(m,e)$ then $\bar{S}=\langle \{m\}+\{x\mod m \mid x\in msg(S) \}\rangle\in\mathcal{C}(m,e)$ and $F(\bar{S})\leq F(S)$.
\end{proposition}

As a consequence of the previous proposition we get the following result.

\begin{corollary}\label{c17}
If $m$ and $e$ are integers such that $m\geq e\geq 2$ then $F(m,e)=\min\{F(S) \mid S\in\mathcal{C}(m,e)\}$.
\end{corollary}

The set $\mathcal{C}(m,e)$ is finite, so previous corollary give us an algorithmic method for computing $F(m,e)$.

\begin{example}
We compute $F(6,5)$. First, we calculate $\mathcal{C}(6,5)$ by using Proposition \ref{p17} and then we apply Corollary \ref{c17}. So,
\begin{multline*}
\mathcal{C}(6,5)=\{\langle 6,7,8,9,10 \rangle, \langle 6,7,8,9,11 \rangle,\\ \langle 6,7,8,10,11 \rangle, \langle 6,7,9,10,11 \rangle, \langle 6,8,9,10,11 \rangle\}
\end{multline*}
and therefore $F(6,5)=\min\{F(\langle 6,7,8,9,10 \rangle)=11$, $F(\langle 6,7,8,9,11 \rangle)=10$, $F(\langle 6,7,8,10,11 \rangle)=9$, $F(\langle 6,7,9,10,11 \rangle)=8$, $F(\langle 6,8,9,10,11 \rangle)=13\}=8$.
\end{example}

Now, we are interested in giving an algorithmic method for computing $\{S\in\mathcal{L}(m,e) \mid F(S)=F(m,e) \}$. Next example shows us that there exist semigroups $S\in\mathcal{L}(m,e)$ such that $S\not\in\mathcal{C}(m,e)$ and $F(S)=F(m,e)$.

\begin{example}The numerical semigroups $S_1=\langle 7,9,10,15 \rangle$ and  $S_2=\langle 7,8,10,19
\rangle$ verify that $S_1,S_2\in \mathcal{L}(7,4)\setminus
\mathcal{C}(7,4)$ and $F(S_1)=F(S_2)=13=F(7,4).$

\end{example}

If $S\in\mathcal{L}(m,e)$ we denote by $\theta(S)$ the numerical semigroup generated by $\{m\}+\{x\mod m \mid x\in msg(S)\}$. Clearly, $\theta(S)\in\mathcal{C}(m,e)$.\\

We define in $\mathcal{L}(m,e)$ the following equivalence relation: $S\mathcal{R}T$ if and only if $\theta(S)=\theta(T)$. We denote by $[S]$ the set $\{T\in\mathcal{L}(m,e) \mid S\mathcal{R}T\}$.  Therefore, the quotient set ${\mathcal{L}(m,e)}/{\mathcal{R}}=\{[S]\mid S\in\mathcal{L}(m,e)\}$ is a partition of $\mathcal{L}(m,e)$. Next result is Theorem 3 from \cite{cadiz}.

\begin{proposition}
Let $m$ and $e$ be integers such that $m\geq e\geq 2$. Then, $\{[S] \mid S\in\mathcal{C}(m,e)\}$ is a partition of $\mathcal{L}(m,e)$. Moreover, if $\{S,T\}\subseteq\mathcal{C}(m,e)$ and $S\neq T$ then $[S]\cap[T]=\emptyset$.
\end{proposition}

As a consequence of the previous proposition we have that for computing $\{S\in\mathcal{L}(m,e) \mid F(S)=F(m,e) \}$ is enough doing the following two steps:
\begin{enumerate}
    \item Compute $A=\{S\in\mathcal{C}(m,e) \mid F(S)=F(m,e) \}$.
    \item For every $S\in A$, compute $\{T\in[S] \mid F(T)=F(S) \}$.
\end{enumerate}
We already know how to compute 1. Now, we focus on giving an algorithm that allows us to compute 2.

Using Algorithm 12 from \cite{cadiz}, for $S\in\mathcal{C}(m,e)$ and $F\in\N$ we get the set $\{T\in[S] \mid F(T)\leq F \}$. Clearly if $S\in\mathcal{C}(m,e)$ then $\{T\in[S] \mid F(T)=F(S)\}=\{T\in [S] \mid F(T)\leq F(S) \}$. Therefore we are going to adapt Algorithm 12 from \cite{cadiz} to our needs for computing 2. But before of that, we introduce some concepts and results.\\

If $S$ is a numerical semigroup, we denote by $M(S)=\max(msg(S))$.
If $S\in\mathcal{C}(m,e)$ we define the graph $G([S])$ as follows: $[S]$ is its set of vertices and $(A,B)\in [S]\times[S]$ is an edge if $msg(B)=(msg(A)\setminus\{M(A)\})\cup\{M(A)-m\}$.

The following result is Theorem 9 from \cite{cadiz}.

\begin{proposition}
If $S\in\mathcal{C}(m,e)$ then $G([S])$ is a tree with root $S$. Moreover, if $P\in [S]$ and $msg(P)=\{n_1<n_2<\cdots<n_e\}$ then the sons of $P$ in $G([S])$ are the numerical semigroups of the form $\langle (\{n_1,\ldots,n_e\}\backslash\{n_k\})\cup\{n_k+n_1\}\rangle$ such that $k\in\{2,\ldots,e\}$, $n_k+n_1>n_e$ and $n_k+n_1\notin\langle \{n_1,\ldots,n_e\}\backslash\{n_k\} \rangle$.
\end{proposition}

Now we give an algorithm such that for a semigroup $S\in\mathcal{C}(m,e)$ it computes the set $\{T\in[S] \mid F(T)=F(S) \}$.

\begin{algorithm}\caption{Sketch of the algorithm to compute the semigroups of each equivalence class such that their Frobenius number is minimum}\label{a32}
\textbf{INPUT:} $S\in\mathcal{C}(m,e)$.\\
\textbf{OUTPUT:} $\{T\in [S] \mid F(T)=F(S) \}$.
\begin{algorithmic}[1]
	%\State $A=\{\langle m,m+1,\ldots,2m-1 \rangle\}$, $I=\emptyset$ and $\alpha=\ceil{\frac{m-1}{e-1}}m-1$.
    \State $A=\{S\}$ and $B=\{S\}$
	\While{True}
		\State $C=\{H \mid H \mbox{ is son of an element of } B \mbox{ and } F(H)=F(S) \}$
		\If{$C=\emptyset$}
			\State \Return A
		\EndIf
		\State $A=A\cup C$, $B=C$.
	\EndWhile
\end{algorithmic}
%    \begin{enumerate}[Step 1.]
%        \item $A=\{S\}$ and $B=\{S\}$.
%        \item $C=\{H \mid H \mbox{ is son of an element of } B \mbox{ and } F(H)=F(S) \}$.
%        \item If $C=\emptyset$ return A.
%        \item $A=A\cup C$, $B=C$ and go to Step 2.
%    \end{enumerate}
\end{algorithm}

We finish this section with an example for illustrate the above algorithm.

\begin{example}
We use now Algorithm \ref{a32} for computing $\{T\in [S] \mid F(T)=F(S)=10 \}$ where $S=\langle 6,7,8,9,11 \rangle\in\mathcal{C}(6,5)$.
    \begin{itemize}
        \item $A=\{\langle 6,7,8,9,11 \rangle\}$ and $B=\{\langle 6,7,8,9,11 \rangle\}$.
        \item $C=\{\langle 6,8,9,11,13 \rangle, \langle 6,8,11,13,15 \rangle\}$.
        \item $A=\{\langle 6,7,8,9,11 \rangle,\langle 6,8,9,11,13 \rangle, \langle 6,8,11,13,15 \rangle\}$ and\\ $B=\{\langle 6,8,9,11,13 \rangle, \langle 6,8,11,13,15 \rangle \}$.
        \item $C=\emptyset$.
    \end{itemize}
Therefore, $\{T\in [S] \mid F(T)=10 \}=\{\langle 6,7,8,9,11 \rangle, \langle 6,8,9,11,13 \rangle, \langle 6,8,11,13,15 \rangle\}$. This result is also obtained by executing the command \texttt{Algorithm3[\{6,7,8,9,11\}]} of \cite{PROGRAMA}).
%Using the package \texttt{FrobeniusNumberAndGenus}, executing the command \texttt{Algorithm3[{6,8,9,11,13}]} we obtain this result.
\end{example}

\end{document}